\documentclass{jaspretmp}

\usepackage{bm}
\usepackage{rotating}
\newcommand{\RR}{\mathbb{R}}
\newcommand{\ZZ}{\mathbb{Z}}
\newcommand{\NN}{\mathbb{N}}

\newcommand{\cE}{{\mathcal E}}
\newcommand{\cF}{{\mathcal F}}

\newcommand{\cP}{{\mathcal P}}

\newcommand{\cU}{{\mathcal U}}

\newcommand{\tri}[3]{\langle #1, #2 \, | \, #3 \rangle}

\newcommand{\bra}[1]{\{ #1 \}}

\newcommand{\indep}{\mathop{\perp\!\!\!\perp}}

\newcommand{\cfg}{\mathcal{U}}
\newcommand{\repr}{\theta}

\newcommand{\tABC}{{\tri{A}{B}{C}}}
\newcommand{\cEABC}{\cE_{\tri{A}{B}{C}}}
\newcommand{\tabG}{{\tri{a}{b}{\Gamma}}}

\newcommand{\abs}[1]{|#1|}

\newcommand{\srift}{slash rift}
\newcommand{\brift}{backslash rift}
\newcommand{\rift}{\mathrm{R}}
\newcommand{\cQ}{{\mathcal Q}}
\newcommand{\cR}{{\mathcal R}}
\newcommand{\GammaSetRepr}{q_{\downarrow}^{\repr}} 
\newcommand{\aGammaSetRepr}{q_{\leftarrow}^{\repr}} 
\newcommand{\bGammaSetRepr}{q_{\rightarrow}^{\repr}} 
\newcommand{\abGammaSetRepr}{q_{\uparrow}^{\repr}} 
\newcommand{\GammaSet}{q_{\downarrow}}
\newcommand{\aGammaSet}{q_{\leftarrow}}
\newcommand{\bGammaSet}{q_{\rightarrow}}
\newcommand{\abGammaSet}{q_{\uparrow}}
\newcommand{\cardA}{\tilde{s}}
\newcommand{\cardB}{\tilde{t}}

\newcommand{\GammaSetSei}{\GammaSet} 
\newcommand{\aGammaSetSei}{\aGammaSet} 
\newcommand{\bGammaSetSei}{\bGammaSet} 

\begin{document}

\title{Properties of semi-elementary imsets as sums of elementary imsets}

\author[T.\ Kashimura, T.\ Sei, A.\ Takemura, K.\ Tanaka]{Takuya Kashimura\affil{1},
Tomonari Sei\affil{2}\affil{3},
Akimichi Takemura\affil{1}\affil{3}\comma\corrauth, 
Kentaro Tanaka\affil{4}}

\address{\affilnum{1}\ Department of Mathematical Informatics,
Graduate School of Information Science and Technology, University of Tokyo,
Tokyo, Japan\\
\affilnum{2}\ Department of Mathematics, Keio University, Tokyo, Japan\\
\affilnum{3}\ JST CREST\\
\affilnum{4}\ Department of Industrial Engineering and Management, Tokyo Institute of Technology, Tokyo, Japan}
\emails{
{\tt kashimura@stat.t.u-tokyo.ac.jp}\ (T.Kashimura),
{\tt takemura@stat.t.u-tokyo.ac.jp}\ (A.Takemura),
{\tt sei@math.keio.ac.jp}\ (T.Sei),
{\tt tanaka.k.al@m.titech.ac.jp}\ (K.Tanaka)
}

\begin{abstract}
We study properties of semi-elementary imsets and elementary imsets
introduced by Studen\'y \cite{stu2005}.  
The rules of the semi-graphoid axiom (decomposition, weak union and contraction) 
for conditional independence statements 
can be translated into a simple identity among three semi-elementary imsets.
By recursively applying the identity, 
any semi-elementary imset can be
written as a sum of elementary imsets, which we call a representation
of the semi-elementary imset. 
A semi-elementary imset has many representations.  
We study properties of the set of possible
representations of a semi-elementary imset
and prove that all representations are connected by 
relations among four elementary imsets.
\end{abstract}



\maketitle

\section{Introduction}
The method of imsets introduced by Studen\'y \cite{stu2005} provides
a very powerful algebraic method for studying 
conditional independence statements which hold under a probability
distribution.  
In this paper we prove some facts 
on semi-elementary imsets when they are represented as non-negative
integer combinations of elementary imsets.  
In particular we prove that all representations of a semi-elementary imset
are connected by 
relations among four elementary imsets.

Let $N$ denote a finite set of random variables and let 
$A,B,C$ denote disjoint subsets of $N$.  
The union $A\cup B$ of two sets $A,B$ is abbreviated as $AB$.
As usual $A\indep B | C$ denotes that the random variables in $A$ are conditionally
independent of those in $B$ given the variables in $C$.  
Three rules of the semi-graphoid axiom, 
i) decomposition, ii) weak union and iii) contraction, can be summarized in the
following single equivalence
\begin{equation}
\label{eq:basic1}
X\indep Y_1Y_2 | Z \ \ \Leftrightarrow\ \ 
X\indep Y_1 | Z    \ \text{and}\  X\indep Y_2 | Y_1Z
\end{equation}
for any disjoint subsets $X,Y_1,Y_2,Z$ of $N$.

\begin{figure}[htbp]
\begin{center}
\includegraphics[scale=1]{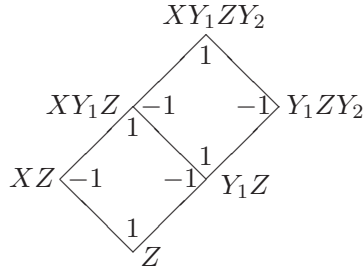}
\end{center}
 \caption{Sum of two semi-elementary imsets}
\label{fig:two-quadrangles}
\end{figure}

Let $\cP(N)$ denote the power set of $N$.
For each triplet of disjoint subsets of $N$, denoted by $\tABC$, 
Studen\'y \cite{stu2005} defined the {\em semi-elementary imset}
$u_{\tABC}:\cP(N)\mapsto \ZZ$ by
\begin{equation}
\label{eq:semi-elementary}
u_{\tABC}(S)=
\begin{cases}
1 & \text{ if } S=ABC \text{ or } S=C \\
-1 & \text{ if } S=AC \text{ or } S=B C \\
0 & \text{otherwise}.
\end{cases}
\end{equation}
In terms of semi-elementary imsets, \eqref{eq:basic1} is written as an identity 
\begin{equation}
\label{eq:basic-u}
u_{\tri{X}{Y_1Y_2}{Z}} = u_{\tri{X}{Y_1}{Z}} + u_{\tri{X}{Y_2}{Y_1Z}} 
\end{equation}
among three semi-elementary imsets. 
From the definition \eqref{eq:semi-elementary} it is easily seen that
\eqref{eq:basic-u} holds.  See Figure \ref{fig:two-quadrangles}.
\eqref{eq:basic-u} is very convenient, because it 
can be regarded as an identity 
among three  $2^{|N|}$-dimensional integer vectors.

When $A$ and $B$ are singletons, written with lower-case letters as $a$ and $b$,
$u_{\tri{a}{b}{C}}$ is called an {\em elementary} imset.
By recursively applying
\eqref{eq:basic-u}, i.e., splitting a semi-elementary imset into two semi-elementary imsets, 
any semi-elementary imset $u_{\tABC}$
can be written as a sum (more precisely a non-negative integer
combination) of elementary imsets 
$u_\tabG$, where $a\in A$, $b\in B$ and $C \subseteq \Gamma\subseteq ABC$.
This corresponds to the fact that elementary imsets
are extreme rays of the cone generated by semi-elementary imsets (Studen\'y \cite{stu2005}).

We call a sum of elementary imsets which is equal to $u_{\tABC}$
a {\em representation} of $u_{\tABC}$ by elementary imsets. 
Depending on the order of applications of \eqref{eq:basic-u}
to various subsets of $A$ and $B$, 
there are many representations of a semi-elementary imset $u_{\tABC}$. 
Furthermore, as we show in this paper, there are representations
which can not be obtained by recursive application of \eqref{eq:basic-u}. 
See Figure~\ref{fig:counterexample} below. 
For a given $\tABC$, let $\cF_{\tABC}$ denote the
set of all possible representations of $u_{\tABC}$.
We call $\cF_{\tABC}$ the $u_{\tABC}$-{\em fiber}.
In this paper we study the properties of this fiber.

Let $X,Y_1,Y_2$ be singletons $a,b_1,b_2$ in \eqref{eq:basic-u} and
let $C \subseteq Z=\Gamma \subseteq ABC$.
Then by changing the roles of $b_1$ and $b_2$, we obtain the following
relation
\begin{equation}
\label{eq:2by2-relation}
u_{\tri{a}{b_1}{\Gamma}} + u_{\tri{a}{b_2}{b_1\Gamma}}=
u_{\tri{a}{b_2}{\Gamma}} + u_{\tri{a}{b_1}{b_2\Gamma}}.
\end{equation}
We call this relation and the similar relation
$
u_{\tri{a_1}{b}{\Gamma}} + u_{\tri{a_2}{b}{a_1\Gamma}}=
u_{\tri{a_2}{b}{\Gamma}} + u_{\tri{a_1}{b}{a_2\Gamma}}
$
a {\em two-by-two basic relation}.
In this paper by a relation we mean an equality between two non-negative
combinations of elementary imsets.
We call the difference of two sides of a two-by-two basic relation
a {\em two-by-two move}. 

We say that two representations $\repr, \repr'\in \cF_{\tABC}$ are
adjacent if they just differ by a two-by-two basic relation, i.e.\ 
$\repr - \repr'$ is a two-by-two 
move.
Furthermore we say that $\cF_{\tABC}$ is {\it connected}
by two-by-two basic relations
if for any two representations $\repr, \repr' \in \cF_{\tABC}$ 
there exists a sequence of representations $\repr=\repr_0, \repr_1, \dots, \repr_K=\repr'$, such that $\repr_{k-1}$ and $\repr_k$ are adjacent, $1\le k\le K$. 
Now our main result is stated as follows.
\begin{theorem}
\label{thm:main}
Every $\cF_{\tABC}$ is connected by two-by-two basic relations.
\end{theorem}

From the viewpoint of toric ideals and Markov bases (e.g.\cite{sturmfels1996},\cite{diaconis-sturmfels},\cite{drton-sturmfels-sullivant-lecture}), 
this result is
closely related to connectivity of a specific fiber by a subset of a Markov basis. 
See \cite{hara-takemura-yoahisda-2010},\cite{chen-dinwoodie-yoshida-2010},%
\cite{yoshida-jalgstat2010} for relevant results.
Since a Markov basis for the whole configuration of elementary imsets
is very complicated  (\cite{hemmecke-etal-2008}), it is
remarkable that $\cF_{\tABC}$ is connected by
the two-by-two basic relations.

The organization of the paper is as follows.
In Section \ref{sec:preliminaries} we set up our notation and 
summarize basic facts on imsets.
In Section \ref{sec:main results} we state our results, including
a sketch of the proof of Theorem \ref{thm:main}.
In Section \ref{sec:computational results}
we show some numerical and computational results on the fiber $\cF_{\tABC}$, when
the cardinalities of $A$ and $B$ are small.
Long proofs of some lemmas and Theorem \ref{thm:main} are given in Section  \ref{sec:proofs}.

\section{Preliminaries}
\label{sec:preliminaries}
In this section we set up our notation and definitions
following \cite{stu2005} and \cite{kashimura-takemura-2011}.
Let $N$ be a finite set  and let $\mathcal{P}(N) = \bra{A \mid A \subseteq N}$ 
denote its power set.
An integer-valued multiset $f:\mathcal{P}(N)\to \ZZ$ is called an imset.
We write the union $A \cup B$ as $AB$.
A singleton set $\{a\}$ is simply written as $a$.

For a triplet $\tABC$
of disjoint subsets of $N$, the
semi-elementary imset $u_{\tABC}$ is defined as 
\eqref{eq:semi-elementary}.  
When $A=a$ and $B=b$ are singletons, $u_{\tri{a}{b}{C}}$ is called elementary.
The set of all elementary imsets for $N$ is denoted as $\cE(N)$.
If $A=\emptyset$ or $B=\emptyset$, then $u_{\tABC}$ is the zero imset.
Hence we usually assume that $A,B$ are non-empty. 
On the other hand $C$ may well be an empty set. 

%

For a given triplet $\tABC$, 
we consider the following set of elementary imsets:
\[
\cEABC = \{u_\tabG \mid a\in A, b\in B, 
C\subseteq \Gamma\subseteq ABC\}.
\]
The cardinality of  $\cEABC$ is given by $|\cEABC|=|A|2^{|A|-1}
|B|2^{|B|-1}$. 
Starting from a given $u_{\tABC}$, consider recursively applying
\eqref{eq:basic-u}.  Then $u_{\tABC}$ is written as
a non-negative integer combination of elementary imsets from $\cEABC$:
\begin{equation}
\label{eq:ABC}
u_{\tABC} = \sum_{u\in \cEABC} \repr_u u, \quad \repr_u\in \NN = \{0,1,2,\dots \}. 
\end{equation}
We call the right-hand side a representation of $u_{\tABC}$.
There are many representations.  
We call the set of possible representations the $u_{\tABC}$-fiber 
and denote it by $\cF_{\tABC}$.
Even for the case $|A|=1$, $|B|=2$,
the two-by-two basic relation in 
\eqref{eq:2by2-relation} shows that there are two representations
of $u_{\tri{a}{b_1 b_2}{C}}$.
If $A=a$ is a singleton, it is easily seen that
there are $|B|!$ different representations of  $u_{\tri{a}{B}{C}}$.  However
for the general  case $|A|\ge 2$, $|B|\ge 2$, it is not trivial to
enumerate $\cF_{\tABC}$.

Now consider writing $u_{\tABC}$ as a $2^{|N|}$-dimensional integer column
vector and $u\in \cEABC$ into a
$2^{|N|}\times |\cEABC|$ integer matrix  $\cU_{\tri{A}{B}{C}}$.
For example 
$\cfg_{\tri{a_1 a_2}{b_1 b_2}{C}}$ is
written as in Table \ref{tab:configuration}.   We call 
$\cfg_{\tABC}$ the configuration for the semi-elementary imset 
$u_{\tABC}$.
\newcommand{\bs}{\begin{sideways}$}
\newcommand{\es}{$\end{sideways}}
\begin{table}[thbp]
\caption{Configuration $\cfg_{\tri{a_1 a_2}{b_1b_2}{C}}$ of elementary imsets in
$\cE_{\tri{a_1 a_2 }{b_1 b_2}{C}}$}
\label{tab:configuration}
{\footnotesize
\setlength{\tabcolsep}{3pt}
\begin{center}
\begin{tabular}{c|cccccccccccccccc}
&
\bs \tri{a_1}{b_1}{a_2 b_2 C} \es &
\bs \tri{a_1}{b_2}{a_2 b_1 C} \es &
\bs \tri{a_2}{b_1}{a_1 b_2 C} \es &
\bs \tri{a_2}{b_2}{a_1 b_1 C} \es &
\bs \tri{a_2}{b_1}{b_2 C} \es &
\bs \tri{a_1}{b_1}{b_2 C} \es &
\bs \tri{a_2}{b_2}{b_1 C} \es &
\bs \tri{a_1}{b_2}{b_1 C} \es &
\bs \tri{a_1}{b_2}{a_2 C} \es &
\bs \tri{a_1}{b_1}{a_2 C} \es &
\bs \tri{a_2}{b_2}{a_1 C} \es &
\bs \tri{a_2}{b_1}{a_1 C} \es &
\bs \tri{a_2}{b_2}{C} \es &
\bs \tri{a_2}{b_1}{C} \es &
\bs \tri{a_1}{b_2}{C} \es &
\bs \tri{a_1}{b_1}{C} \es 
\\ \hline
 $a_1 a_2 b_1 b_2C$ &
                    1& 1& 1& 1& 0& 0& 0& 0& 0& 0& 0& 0& 0& 0& 0& 0 \\
 $a_2 b_1 b_2C$ &              	      	    	  	  
                    -1& -1& 0& 0& 1& 0& 1& 0& 0& 0& 0& 0& 0& 0& 0& 0 \\
 $a_1 b_1 b_2C$ &              	      	    	  	  
                    0& 0& -1& -1& 0& 1& 0& 1& 0& 0& 0& 0& 0& 0& 0& 0 \\
 $a_1 a_2 b_2C$ &              	      	    	  	  
                    -1& 0& -1& 0& 0& 0& 0& 0& 1& 0& 1& 0& 0& 0& 0& 0 \\
 $a_1 a_2 b_1C$ &              	      	    	  	  
                    0& -1& 0& -1& 0& 0& 0& 0& 0& 1& 0& 1& 0& 0& 0& 0 \\
 $b_1 b_2C$ &                	      	    	  	  
                    0& 0& 0& 0& -1& -1& -1& -1& 0& 0& 0& 0& 0& 0& 0& 0 \\
 $a_2 b_2C$ &                	      	    	  	  
                    1& 0& 0& 0& -1& 0& 0& 0& -1& 0& 0& 0& 1& 0& 0& 0 \\
 $a_2 b_1C$ &                	      	    	  	  
                    0& 1& 0& 0& 0& 0& -1& 0& 0& -1& 0& 0& 0& 1& 0& 0 \\
 $a_1 b_2C$ &                	      	    	  	  
                    0& 0& 1& 0& 0& -1& 0& 0& 0& 0& -1& 0& 0& 0& 1& 0 \\
 $a_1 b_1C$ &                	      	    	  	  
                    0& 0& 0& 1& 0& 0& 0& -1& 0& 0& 0& -1& 0& 0& 0& 1 \\
 $a_1 a_2C$ &                	      	    	  	  
                    0& 0& 0& 0& 0& 0& 0& 0& -1& -1& -1& -1& 0& 0& 0& 0 \\
 $b_2C$ &                    	      	    	  	  
                    0& 0& 0& 0& 1& 1& 0& 0& 0& 0& 0& 0& -1& 0& -1& 0 \\
 $b_1C$ &                    	      	    	  	  
                    0& 0& 0& 0& 0& 0& 1& 1& 0& 0& 0& 0& 0& -1& 0& -1 \\
 $a_2C$ &                    	      	    	  	  
                    0& 0& 0& 0& 0& 0& 0& 0& 1& 1& 0& 0& -1& -1& 0& 0 \\
 $a_1C$ &                    	      	    	  	  
                    0& 0& 0& 0& 0& 0& 0& 0& 0& 0& 1& 1& 0& 0& -1& -1 \\
 $C$& 		             	      	    	  	  
                    0& 0& 0& 0& 0& 0& 0& 0& 0& 0& 0& 0& 1& 1& 1& 1 \\
\end{tabular}
\end{center}
}
\end{table}
Then a representation \eqref{eq:ABC} is written in a matrix form as
\[
u_{\tABC} = \cfg_{\tABC} \repr,
\]
where $\repr$ is the column vector of coefficients $\repr_u$ in \eqref{eq:ABC}.
From now on we identify  a non-negative integer combination of elementary
imsets from $\cEABC$ with the vector of non-negative
integer coefficients.  Then $\cF_{\tABC}$ is written as
\begin{equation}
\label{eq:fiber-vector}
\cF_{\tABC} =
\{ \repr\in \NN^{|\cEABC|} \mid
u_{\tABC} = \cfg_{\tABC} \repr\}.
\end{equation}
In this form it is evident that $\cF_{\tABC}$ is a particular fiber
in the theory of Markov basis.

Note that the configuration $\cfg_{\tABC}$ is a
subconfiguration of the set of all elementary imsets $\cE(N)$.
It is known that the elementary imsets in $\cEABC$ are
the extreme rays of a face of the cone generated by all the elementary imsets. 
Hence the subconfiguration $\cfg_{\tABC}$ generates a combinatorial 
pure subring in the sense of Ohsugi, Herzog and Hibi (\cite{ 
ohsugi-herzog-hibi2000,ohsugi-2007,ohsugi-hibi2010}).

From the form of $\cfg_{\tABC}$ in Table \ref{tab:configuration},
it is evident that the structure of the 
fiber $\cF_{\tABC}$ only depends on $|A|$ and $|B|$.  
Therefore for our study of the structure of the fiber
$\cF_{\tABC}$, we can assume that $C=\emptyset$ and $N=AB$, without
loss of generality.  

\begin{figure}
\begin{center}
\includegraphics[scale=1]{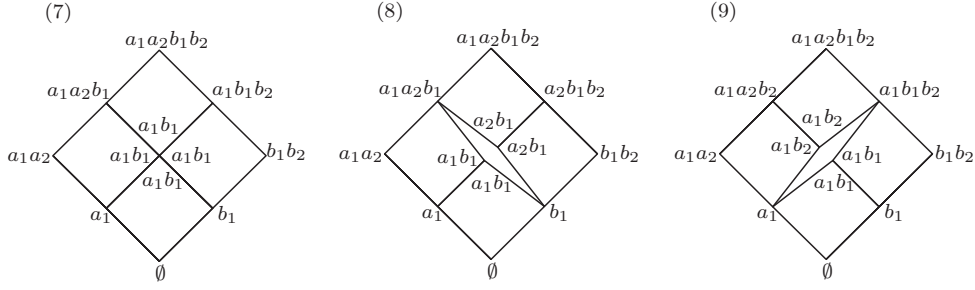}
\end{center}
\caption{Three types of representations of $u_{\tri{a_1 a_2}{b_1 b_2}{\emptyset}}$}
\label{fig:a1a2b1b2}
\end{figure}

For the rest of this section, for the purpose of illustration, we write out
the fiber $\cF_{\tri{a_1 a_2 }{b_1 b_2 }{\emptyset}}$.  If we split $b_1 b_2$ first, we have
$u_{\tri{a_1 a_2}{b_1 b_2}{\emptyset}}=u_{\tri{a_1 a_2}{b_1}{\emptyset}} + u_{\tri{a_1 a_2}{b_2}{b_1}}$
or
$u_{\tri{a_1 a_2}{b_1 b_2}{\emptyset}}=u_{\tri{a_1 a_2}{b_2}{\emptyset}} + u_{\tri{a_1 a_2}{b_1}{b_2}}$.
Consider the former case.  Then we can split $a_1 a_2$, independently in two terms on 
the right-hand side.  Then we have four representations:
\begin{eqnarray}
u_{\tri{a_1 a_2}{b_1 b_2}{\emptyset}}
&=u_{\tri{a_1}{b_1}{\emptyset}} + u_{\tri{a_2}{b_1}{a_1}}
+ u_{\tri{a_1}{b_2}{b_1}} +  u_{\tri{a_2}{b_2}{a_1 b_1}}
 \label{eq:a2b2_type1}
\\
&=u_{\tri{a_2}{b_1}{\emptyset}} + u_{\tri{a_1}{b_1}{a_2}}
 + u_{\tri{a_1}{b_2}{b_1}} +  u_{\tri{a_2}{b_2}{a_1 b_1}}
\nonumber 
\\
&=
u_{\tri{a_1}{b_1}{\emptyset}} + u_{\tri{a_2}{b_1}{a_1}}
+ u_{\tri{a_2}{b_2}{b_1}} +  u_{\tri{a_1}{b_2}{a_2 b_1}}
\label{eq:a2b2_type2}
\\
&=
u_{\tri{a_2}{b_1}{\emptyset}} + u_{\tri{a_1}{b_1}{a_2}}
+ u_{\tri{a_2}{b_2}{b_1}} +  u_{\tri{a_1}{b_2}{a_2 b_1}}.
\nonumber 
\end{eqnarray}
When we enumerate the representations, apparently we have 16 representations.  However
in fact there are 12 distinct representations of $\cF_{\tri{a_1 a_2 }{b_1 b_2 }{\emptyset}}$.
If we consider symmetry with respect to the interchanges $a_1 \leftrightarrow a_2$ and
$b_1 \leftrightarrow b_2$, there are three types of representations as shown in 
Figure \ref{fig:a1a2b1b2}.  Note that two types of representations have ``rifts'', which will
be discussed in Sections \ref{sec:main results} and \ref{sec:proofs}. 
The representations of the first and second types in Figure~\ref{fig:a1a2b1b2} are 
obtained in (\ref{eq:a2b2_type1}) and (\ref{eq:a2b2_type2}). 
By splitting $a_{1}a_{2}$ first, 
the representation of the third type in Figure~\ref{fig:a1a2b1b2} is given as follows:
\begin{equation}
 u_{\tri{a_1 a_2}{b_1 b_2}{\emptyset}}
  =
  u_{\tri{a_1}{b_1}{\emptyset}} + u_{\tri{a_2}{b_2}{a_1}}
  + u_{\tri{a_1}{b_2}{b_1}} +  u_{\tri{a_2}{b_1}{a_1 b_2}}. 
  \label{eq:a2b2_type3} 
\end{equation}
The type without rifts
is obtained by applying the same splitting to two intermediate terms
on the right-hand side of 
$u_{\tri{a_1 a_2}{b_1 b_2}{\emptyset}}=u_{\tri{a_1 a_2}{b_1}{\emptyset}} + u_{\tri{a_1 a_2}{b_2}{b_1}}$
or
$u_{\tri{a_1 a_2}{b_1 b_2}{\emptyset}}=u_{\tri{a_1 a_2}{b_2}{\emptyset}} + u_{\tri{a_1 a_2}{b_1}{b_2}}$.
This is an important fact in proving Theorem \ref{thm:main}.

\section{Main results}
\label{sec:main results}

In this section we establish some facts on the fiber of representations
$\cF_{\tri{A}{B}{C}}$.  At the end of the section we give a sketch of a proof of
Theorem \ref{thm:main}.  Some long proofs are given in Section
\ref{sec:proofs}.

First we show that a representation always involves $|A|\times |B|$ elementary imsets.
This fact is obvious if a representation is obtained by recursive application of
\eqref{eq:basic-u}.  However since there are representations not obtained 
by recursive application of
\eqref{eq:basic-u}, we need to show
this fact for all representations.

As before we consider an elementary imset $u\in \cEABC$ as
a function from $\cP(N)$ to $\ZZ \subset \RR$ or a $2^{|N|}$-dimensional
integer vector.  
We define the inner product of two functions $f,g$ from 
$\cP(N)$ to $\RR$ by the standard inner product of $2^{|N|}$-dimensional
vectors.
\[
\langle f, g\rangle = \sum_{S\in \cP(N)} f(S)g(S)=f(\emptyset)g(\emptyset) + \cdots + 
f(N)g(N).
\]
%
%

For $u=u_{\tabG}\in \cEABC$, let $s=|A\cap \Gamma|, t=|B\cap \Gamma|$. 
Then $0 \le s \le |A|-1$, $0\le t\le |B|-1$.
Let
\[
\cEABC^{s,t} = \{ u\in \cEABC \mid s=|A\cap \Gamma|, t=|B\cap \Gamma|\}.
\]
Classifying the terms in a representation, we can rewrite \eqref{eq:ABC} as
\[
u_{\tri{A}{B}{C}}=\sum_{s=0}^{|A|-1} \sum_{t=0}^{|B|-1} \sum_{u\in \cEABC^{s,t}} \repr_u u.
\]
%
We now prove that a representation is a sum of $|A||B|$ terms and each 
term is taken just once from $\cEABC^{s,t}$.
\begin{proposition}
\label{prop:|A||B|}
Let $\theta: u_{\tABC}=\sum_{u\in \cEABC} \repr_u u$ be a representation.
Then for each $(s,t)$, $0\le s \le |A|-1$, $0\le t\le |B|-1$, 
there is one $u\in \cEABC^{s,t}$ such that $\repr_u=1$ and 
$\repr_{u'}=0$ for other $u'\in \cEABC^{s,t}$, $u'\neq u$.
\end{proposition}

\begin{proof}
For $0\le s \le |A|-1$ and for $S\subseteq A$ let $g_s: \cP(A) \rightarrow \ZZ$
be defined as
\[
g_s(S)=1_{\{|S|>s\}}= \begin{cases} 1 & \text{ if } |S| > s \\
     0 & \text{ otherwise}.
   \end{cases}
\]
Let $h_t: \cP(B)\rightarrow \ZZ$ be similarly defined. 
Let $f(S)=g_s(A\cap S)h_t(B\cap S)$.  We consider the inner product of $f$
with $u=u_{\tri{a}{b}{\Gamma}} \in \cEABC$.
Write $A'=A\cap \Gamma$, $B'=B\cap \Gamma$. 
It is easily seen that 
the inner product of $f$ with $u_{\tri{a}{b}{\Gamma}}$ is given as
\[
\langle f, u_{\tri{a}{b}{\Gamma}} \rangle =
(1_{\{|aA'|>s\}} - 1_{\{|A'|>s\}})(1_{\{|bB'|>t\}} - 1_{\{|B'|>t\}}).
\]
Hence
\[
\langle f, u \rangle = \begin{cases} 1 & \text{ if }  u\in \cEABC^{s,t}\\
   0 & \text{ otherwise }.
 \end{cases}
\]
On the other hand $\langle u_\tABC , f \rangle = (1-0)(1-0)=1$.  Hence  we have
\[
1 = \sum_{u\in \cEABC^{s,t}} \repr_u.
\]
Since $\repr_u$ are non-negative integers we have the proposition.
\qed
\end{proof}

From this proposition, given a representation $\repr: u_\tABC = \sum  \repr_u u$,
for each $0\le s\le |A|-1$, $0\le t\le |B|-1$, there exists
$u(s,t)=u_\repr(s,t)\in \cEABC^{s,t}$ such that 
\begin{equation}
\label{eq:ust}
u_\tABC = \sum_{s=0}^{|A|-1} \sum_{t=0}^{|B|-1} u_\repr(s,t) .
\end{equation}
We often omit the subscript $\repr$ in $u_\repr(s,t)$ and simply write
$u(s,t)$ for notational simplicity.

We say that a representation (\ref{eq:ust}) is $\sigma$-decomposable
if it is obtained by recursive application of (\ref{eq:basic-u}).
Here $\sigma$ stands for ``semi-graphoid''.
A precise definition of $\sigma$-decomposability is given as follows.
It is analogous to the definition of decomposable graphs.
We first define separability of a representation.

\begin{definition}
 Let $1\le s_0\le |A|-1$.
 Then we say that a representation (\ref{eq:ust}) is $(s_0,*)$-separable if there is some
 $A_0\subseteq A$ such that $|A_0|=s_0$ and
 \begin{equation}
  \label{eq:s-sep}
  u_{\tri{A_0}{B}{C}} = \sum_{s=0}^{s_0-1}\sum_{t=0}^{|B|-1}u(s,t),
   \quad u_{\tri{A\setminus A_0}{B}{A_0C}} = \sum_{s=s_0}^{|A|-1}\sum_{t=0}^{|B|-1}u(s,t).
 \end{equation}
 Similarly, for $1\le t_0\le |B|-1$, we say that a representation (\ref{eq:ust})
 is $(*,t_0)$-separable if there is some $B_0\subseteq B$ such that $|B_0|=t_0$ and
 \begin{equation}
  \label{eq:t-sep}
  u_{\tri{A}{B_0}{C}} = \sum_{s=0}^{|A|-1}\sum_{t=0}^{t_0-1}u(s,t),
   \quad u_{\tri{A}{B\setminus B_0}{B_0C}} = \sum_{s=0}^{|A|-1}\sum_{t=t_0}^{|B|-1}u(s,t).
 \end{equation}
\end{definition}

For example, a representation
\begin{align*}
 u_{\tri{a_1a_2}{b_1b_2}{\emptyset}}
 &= u(0,0) + u(1,0) + u(0,1) + u(1,1)
 \\
 &= u_{\tri{a_1}{b_1}{\emptyset}} + u_{\tri{a_2}{b_1}{a_1}}
 + u_{\tri{a_2}{b_2}{b_1}} + u_{\tri{a_1}{b_2}{a_2b_1}}
\end{align*}
is $(*,1)$-separable with $B_0=b_1$, but not $(1,*)$-separable.

Now $\sigma$-decomposability is recursively defined as follows.

\begin{definition}
 \label{def:sigma_decomposability}
 We say that a representation (\ref{eq:ust}) of $u_{\tri{A}{B}{C}}$ is {\em $\sigma$-decomposable} if
 at least one of the following conditions is satisfied:
 \begin{itemize}
  \item[(i)] $|A|=|B|=1$.
  \item[(ii)] There is some $1\le s_0\le |A|-1$ such that (\ref{eq:ust}) is $(s_0,*)$-separable
	     and two representations in (\ref{eq:s-sep}) are respectively $\sigma$-decomposable.
  \item[(iii)] There is some $1\le t_0\le |B|-1$ such that (\ref{eq:ust}) is $(*,t_0)$-separable
	     and two representations in (\ref{eq:t-sep}) are respectively $\sigma$-decomposable.
 \end{itemize}
\end{definition}
If a representation is not $\sigma$-decomposable, we call the representation {\em $\sigma$-indecomposable}. 
The following proposition will be proved in Section~\ref{sec:proofs}.
\begin{proposition}
\label{prop:sufficient_condition_for_sigma_decomposability}
All of the representations of $u_{\tABC}$  by elementary imsets in $\cEABC$
are $\sigma$-decomposable if and only if  $|A|\le 2$ or $|B|\le 2$.
\end{proposition}



\begin{figure}[htbp]
\begin{center}
\includegraphics[scale=1]{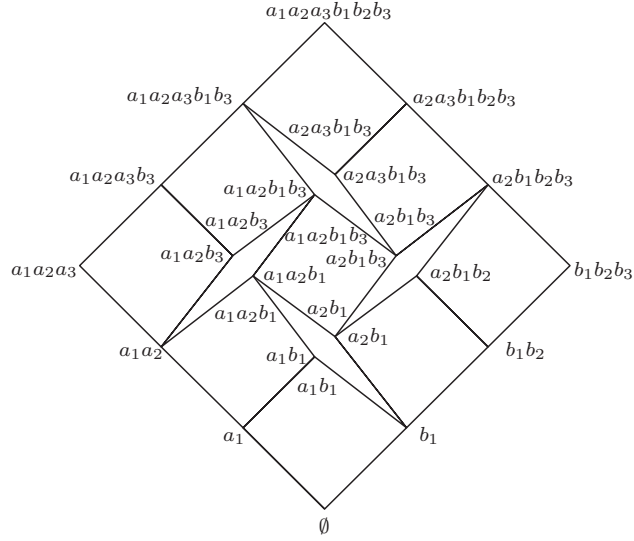}
\end{center}
\caption{A counter-example which is not $\sigma$-decomposable}
\label{fig:counterexample}
\end{figure}

We now consider the case of $u_{\tABC}$ where $\abs{A} \geq 3$ and $\abs{B} \geq 3$. 
In this case, we can construct counter-examples, i.e., representations of $u_{\tABC}$ by elementary imsets which are $\sigma$-indecomposable. 
Now, let us consider the following representation of $u_{\tri{a_{1}a_{2}a_{3}}{b_{1}b_{2}b_{3}}{\emptyset}}$: 
\begin{eqnarray}
 u_{\tri{a_{1}a_{2}a_{3}}{b_{1}b_{2}b_{3}}{\emptyset}}
  & = &
  u_{\tri{a_{1}}{ b_{1}}{\emptyset}}
  +
  u_{\tri{ a_{2}}{ b_{1}}{a_{1}}}
  +
  u_{\tri{ a_{3}}{ b_{3}}{a_{1}a_{2}}}
  \nonumber \\
  & & +
  u_{\tri{ a_{2}}{ b_{2}}{b_{1}}}
  +
  u_{\tri{ a_{1}}{ b_{3}}{a_{2}b_{1}}}
  +
  u_{\tri{ a_{3}}{ b_{1}}{a_{1}a_{2}b_{3}}}
  \nonumber \\
  & & +
  u_{\tri{ a_{2}}{ b_{3}}{b_{1}b_{2}}}
  +
  u_{\tri{ a_{3}}{ b_{2}}{a_{2}b_{1}b_{3}}}
  +
  u_{\tri{ a_{1}}{ b_{2}}{a_{2}a_{3}b_{1}b_{3}}}. 
 \label{eq:counter_example}
 \end{eqnarray}
See Figure~\ref{fig:counterexample}. 
From the definition of $\sigma$-decomposability and Figure~\ref{fig:counterexample}, 
it is clear that the representation of (\ref{eq:counter_example}) is 
$\sigma$-indecomposable. 
Furthermore, from Proposition~\ref{prop:sufficient_condition_for_sigma_decomposability}, 
the counter-example given by~(\ref{eq:counter_example}) is the smallest representation which is $\sigma$-indecomposable. 
The counter-example (\ref{eq:counter_example}) can be extended to
$|A|>3$ or $|B|>3$ as in Figure \ref{fig:counterexample_examples} below.

In the rest of this section, 
we give only a sketch of a proof of Theorem~\ref{thm:main}. 
We will give the complete proof in Section~\ref{sec:proofs}. 

First we note that every representation of $u_{\tABC}$ has its rift pattern.
Intuitively, a rift pattern is a rectangle consisting of $|A|\times |B|$ cells with some rifts,
where each cell corresponds to an elementary imset $u(s,t)$ in \eqref{eq:ust}. 
For example, \eqref{eq:counter_example} has a rift pattern with 
nine cells and four rifts (of length 2) as shown in Figure \ref{fig:counterexample}.
A precise definition of rifts and rift patterns will be given in Definition \ref{def:rift}.

We label elements of $A$ and $B$ as $A=\{a_1,\ldots,a_{|A|}\}$ and $B=\{b_1,\ldots,b_{|B|}\}$.
As discussed in the last paragraph of Section~\ref{sec:preliminaries}, 
representations without rifts are obtained by applying the same splitting rules of~(\ref{eq:basic-u}) to the related intermediate terms recursively. 
In particular,
we can obtain the following representation without rifts:
\begin{equation}
 \label{eq:standard_represetation}
 u_{\tABC} = \sum_{i=1}^{\abs{A}}\sum_{j=1}^{\abs{B}} u_{\tri{a_{i}}{b_{j}}{a_{1} \dots a_{i-1} b_{1} \dots b_{j-1}C}}.
\end{equation}
We call the above representation the {\em standard representation} of $u_{\tABC}$. 
It is easily seen that every other representations without rifts are obtained 
by permuting $a_i$'s and $b_j$'s independently in the standard representation. 
The converse is also true, i.e., 
any representations obtained by permuting $a_i$'s and $b_j$'s independently in the standard representation have no rifts. 
By the following proposition, it suffices to show that we can eliminate all of the rifts in the representations by applying two-by-two basic relations. 
\begin{proposition}
 \label{prop:kashimura}
 Any representations of $u_{\tABC}$ without rifts are obtained by applying two-by-two basic relations to the standard representation in~(\ref{eq:standard_represetation}). 
\end{proposition}
\begin{proof}
 First, note that any permutation is constructed from adjacent transpositions $a_{i} \leftrightarrow a_{i+1}$ and $b_{j} \leftrightarrow b_{j+1}$. 
 Therefore, we only need to show that any adjacent transposition is constructed by two-by-two basic relations. 
 For example, suppose that we want to interchange $a_{1} \leftrightarrow a_{2}$ in the standard representation in~(\ref{eq:standard_represetation}). 
 Then, this adjacent transposition is realized by applying $\abs{B}$ two-by-two basic relations: 
 \begin{eqnarray}
 u_{\tri{a_{1}}{b_{1}}{C}} + u_{\tri{a_2}{b_{1}}{a_{1}C}}
  & = &
  u_{\tri{a_{2}}{b_{1}}{C}} + u_{\tri{a_1}{b_{1}}{a_{2}C}},
  \nonumber \\
 u_{\tri{a_{1}}{b_{2}}{b_{1}C}} + u_{\tri{a_2}{b_{2}}{a_{1}b_{1}C}}
  & = &
  u_{\tri{a_{2}}{b_{2}}{b_{1}C}} + u_{\tri{a_1}{b_{2}}{a_{2}b_{1}C}},
  \nonumber \\
 & \vdots & 
  \nonumber \\ 
 u_{\tri{a_{1}}{b_{\abs{B}}}{b_{1} \dots b_{\abs{B}-1}C}} + u_{\tri{a_2}{b_{\abs{B}}}{a_{1}b_{1} \dots b_{\abs{B}-1}C}}
  & = &
  u_{\tri{a_{2}}{b_{\abs{B}}}{b_{1} \dots b_{\abs{B}-1}C}} + u_{\tri{a_1}{b_{\abs{B}}}{a_{2}b_{1} \dots b_{\abs{B}-1}C}}. 
  \nonumber 
 \end{eqnarray} 
 In the same way, we can show that any adjacent transpositions $a_{i} \leftrightarrow a_{i+1}$ in arbitrary representations without rifts 
 are constructed by using two-by-two basic relations. 
 The proof for any adjacent transpositions $b_{j} \leftrightarrow b_{j+1}$ is also similar to the proof for $a_{i}$'s. 
 \qed
\end{proof}

\begin{figure}
\begin{center}
\includegraphics[scale=1]{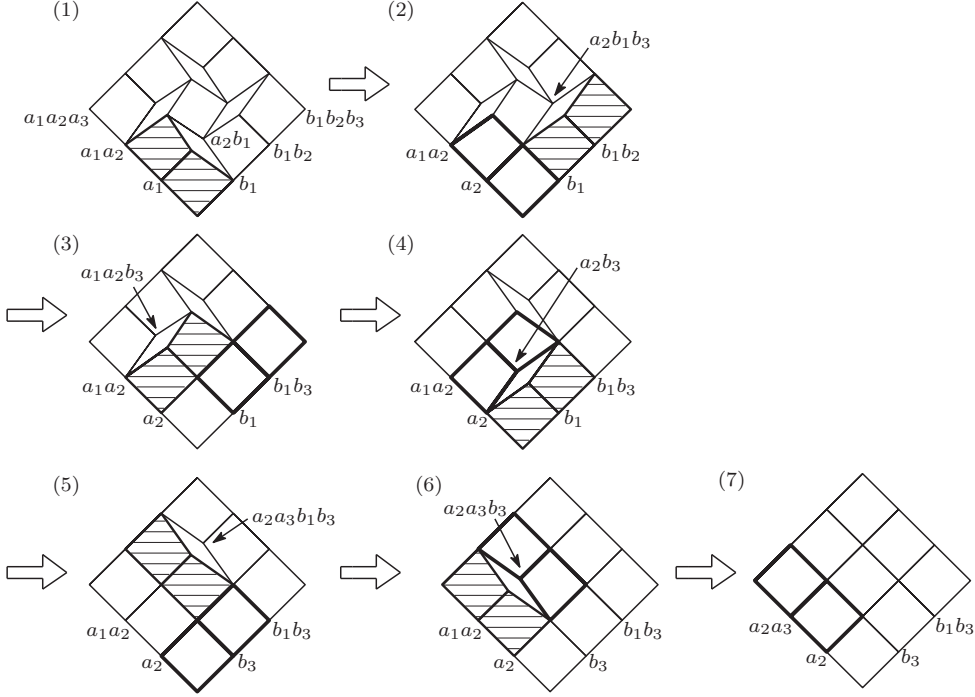}
\end{center}
\caption{Elimination of rifts}
\label{fig:elimination_of_rifts_in_counterexample}
\end{figure}
Rather than stating the method to eliminate the rifts in general, let us work through an example to see how we can do it. 
Here, we use 
the example of (\ref{eq:counter_example})
to explain the method to eliminate the rifts. 
First, by applying a two-by-two basic relation
$u_{\tri{a_{1}}{b_{1}}{\emptyset}} + u_{\tri{a_{2}}{b_{1}}{a_{1}}} = u_{\tri{a_{2}}{b_{1}}{\emptyset}} + u_{\tri{a_{1}}{b_{1}}{a_{2}}}$, 
we can change the representation (1) in  Figure~\ref{fig:elimination_of_rifts_in_counterexample} 
into (2) in Figure~\ref{fig:elimination_of_rifts_in_counterexample} and eliminate the lower left rift. 
Note that the shaded areas in Figure~\ref{fig:elimination_of_rifts_in_counterexample} 
represent the elementary imsets to which we are applying two-by-two basic relations. 
In the same way, we can eliminate the lower right rift by applying
$u_{\tri{a_{2}}{b_{2}}{b_{1}}} + u_{\tri{a_{2}}{b_{3}}{b_{1}b_{2}}} = u_{\tri{a_{2}}{b_{3}}{b_{1}}} + u_{\tri{a_{2}}{b_{2}}{b_{1}b_{3}}}$. 
The third step changing from (3) to (5) in Figure~\ref{fig:elimination_of_rifts_in_counterexample} is rather complicated. 
We can eliminate the upper left rift by applying 
$u_{\tri{a_{1}}{b_{1}}{a_{2}}} + u_{\tri{a_{1}}{b_{3}}{a_{2}b_{1}}} = u_{\tri{a_{1}}{b_{3}}{a_{2}}} + u_{\tri{a_{1}}{b_{1}}{a_{2}b_{3}}}$. 
However, this operation also generates a lower left rift as in (4) of Figure~\ref{fig:elimination_of_rifts_in_counterexample}. 
In this case, we must apply further two-by-two basic relation
$u_{\tri{a_{2}}{b_{1}}{\emptyset}} + u_{\tri{a_{2}}{b_{3}}{b_{1}}} = u_{\tri{a_{2}}{b_{3}}{\emptyset}} + u_{\tri{a_{2}}{b_{1}}{b_{3}}}$, 
to eliminate the rift generated by side effects of a two-by-two basic relation. 
The operation of the fourth step changing from (5) to (7) in Figure~\ref{fig:elimination_of_rifts_in_counterexample} 
is almost the same as in the third step. 
In this manner, we can eliminate all the rifts in the representations. 
Note that, in the first step,  we can not apply 
$u_{\tri{a_{1}}{b_{1}}{a_{2}}} + u_{\tri{a_{1}}{b_{3}}{a_{2}b_{1}}} = u_{\tri{a_{1}}{b_{3}}{a_{2}}} + u_{\tri{a_{1}}{b_{1}}{a_{2}b_{3}}}$ 
which is used in the third step 
to eliminate the upper left rift. 
This means that order of eliminations of the rifts is important. 
The rigorous definition of these operations and the detailed proofs about the properties are given in Section~\ref{sec:proofs}. 

\section{Some numerical and computational results}
\label{sec:computational results}

We report numerical and computational study to 
enumerate all of the representations of $u_{\tABC}$. 
In this section, 
we identify two representations 
that one is obtained by permuting $a_i$'s and $b_j$'s independently in another,
and count the number of representatives of the representations under this equivalence relation. 
Note that if the permutation is not identity, the two representations are different from each other
since they are distinguished by, for example, $u(s,0)$'s and $u(0,t)$'s in \eqref{eq:ust}.
Hence each equivalence class has $\abs{A}!\abs{B}!$ elements.
We call the set of all of the rifts of a representation a {\em rift pattern}.
See Definition \ref{def:rift} for the precise definition. 
Note that two representations may have the same rift pattern. 

First, let us consider the number of rift patterns. 
As will be explained in Section~\ref{sec:proofs}, 
rift occurs at $(\abs{A}-1)(\abs{B}-1)$ points in $\mathcal{P}(N)$
and rifts never intersect one another, i.e., each point in $\mathcal{P}(N)$ has only three patterns: no rift, \srift\ and \brift. 
The precise definitions of \srift\ and \brift\ are given in Section~\ref{sec:proofs} (See Figure~\ref{fig:rift} below). 
Therefore, the number of rift patterns is given by $3^{(\abs{A}-1)(\abs{B}-1)}$. 
The list of all the rift patterns are also given by this enumeration. 

For each rift pattern, we can enumerate all the representations which are consistent with the rift pattern
by making an exhaustive investigation. 
If we only want to count the number of representations, 
then we can compute it as follows. 
First, note that the diversity of representations is given by rifts, 
and, for each rift, the degree of freedom is decided only by the ``length of the rift''. 
Intuitively, in Figure~\ref{fig:rift} in Section~\ref{sec:proofs}, 
the length of the rift is given by $t_{U}-t_{L}$ for a \srift\ 
and $s_{U}-s_{L}$ for a \brift. 
The following lemma is useful for counting the degree of freedom. 
\begin{proposition}
 \label{prop:deg_of_free}
 For any rift with length $l$, we can compute the degree of freedom $d(l)$ 
 by the following recurrence formula: 
 \begin{equation}
  d(l) = l! - \sum_{k=1}^{l-1} (l-k)! \, d(k), \ \ \ l\geq 2,
  \label{eq:deg_of_free}
 \end{equation}
 where $d(1) = 1$.
\end{proposition}
We will give the proof of the above proposition in Section~\ref{sec:proofs}. 
Since the diversity increases in a multiplicative fashion independently for each rift, 
the number of all of the representations for a rift pattern 
is obtained by taking the products of the degree of freedoms for all of the rifts. 
Hence, by summing the number of representations for all of the rift patterns, 
we can count the number of representations. 

The number of $\sigma$-indecomposable representations is calculated by applying the recursive algorithm 
given in Definition~\ref{def:sigma_decomposability} to each rift pattern. 

If $\abs{A}=1$ or $\abs{B}=1$, then there is no rift.
Therefore, the numbers of rift patterns and representations of $u_{\tABC}$ for $\abs{A} = 1$ or $\abs{B} = 1$ are both 1. 
Furthermore, from Proposition~\ref{prop:sufficient_condition_for_sigma_decomposability}, 
the number of $\sigma$-indecomposable representations of $u_{\tABC}$ for $\abs{A} \leq 2$ or $\abs{B} \leq 2$ is zero. 

Next, let us consider the case for $\abs{A}=2$. 
In this case, we have the following lemma.
\begin{proposition}
 \label{prop:num_of_rep_for_a2}
 Let $r_{2}(\abs{B})$ denote the number of representatives of the representations of $u_{\tABC}$ for $\abs{A}=2$. 
 Then, we have the following recurrence formula: 
 \begin{equation}
  r_{2}(m) = m! + \sum_{k=1}^{m-1} k! \, r_{2}(m-k), \ \ \ m\geq 2,
  \label{eq:num_of_rep_for_a2}
 \end{equation}
 where $r_{2}(1) = 1$. 
\end{proposition}
We will give the proof of the above proposition in Section~\ref{sec:proofs}. 

\begin{figure}
\begin{center}
\includegraphics[scale=1]{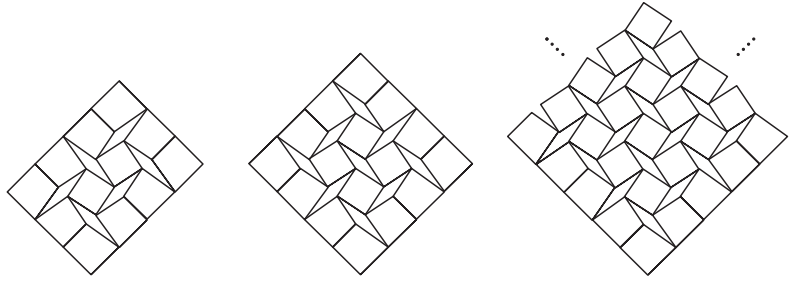}
\end{center}
\caption{$\sigma$-indecomposable rift patterns}
\label{fig:counterexample_examples}
\end{figure}

For the case of $2 \leq \abs{A}, \abs{B} \leq 5$. 
we give the numbers of rift patterns, representations, $\sigma$-indecomposable rift patterns, $\sigma$-indecomposable representations 
of $u_{\tABC}$ 
in Table~\ref{tab:computaional_results}. 
We also give some examples of $\sigma$-indecomposable rift patterns in Figure~\ref{fig:counterexample_examples}.
The left-hand side of Figure~\ref{fig:counterexample_examples} is an example of $\sigma$-indecomposable rift pattern 
for $\abs{A}=3, \abs{B}=4$. 
This can be seen as an extension of the counter-example of Figure~\ref{fig:counterexample} for $\abs{A}=3, \abs{B}=3$. 
The example in the middle of Figure~\ref{fig:counterexample_examples} is an example of $\sigma$-indecomposable rift pattern 
for $\abs{A}=4, \abs{B}=4$
and also can be seen as an extension of the counter-example of Figure~\ref{fig:counterexample}. 
It is easily seen that we can construct an example of $\sigma$-indecomposable rift pattern in general case 
by extending the counter-example of Figure~\ref{fig:counterexample} 
as in the right-hand side of Figure~\ref{fig:counterexample_examples}. 

\begin{table}[thbp]
 \caption{Number of rift patterns, representations, $\sigma$-indecomposable rift patterns and $\sigma$-indecomposable representations for each $\abs{A}$ and $\abs{B}$}
 \label{tab:computaional_results}
 \begin{center}
  \begin{tabular}[c]{cccccc}
   {$\abs{A}$}
   & $\abs{B}$
   & rift patterns
   & \shortstack{representations 
   }
   & \shortstack{$\sigma$-indecomposable \\ rift patterns 
   }
   & \shortstack{$\sigma$-indecomposable \\ representations 
   }
   \\ \hline 
  2 & 2 & 3 & 3 & 0 & 0 \\
  2 & 3 & 9 & 11 & 0 & 0 \\
  2 & 4 & 27 & 47 & 0 & 0 \\
  2 & 5 & 81 & 231 & 0 & 0 \\
  3 & 3 & 81 & 161 & 2 & 2 \\
  3 & 4 & 729 & 2971 & 40 & 96 \\
  3 & 5 & 6561 & 69281 & 562 & 3582 \\
  4 & 4 & 19683 & 241291 & 2436 & 19996 \\
  4 & 5 & 531441 & 25897047 & 102576 & 3420918 \\
  5 & 5 & 43046721 & 12606896129 & 12833474 & 2714509138 \\
  \end{tabular}
 \end{center}
\end{table}

\section{Proofs}
\label{sec:proofs}

\begin{figure}
\begin{center}
\includegraphics[scale=1]{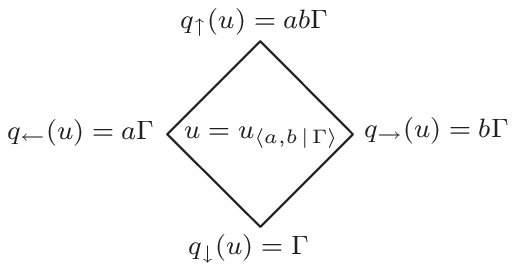}
\end{center}
\caption{Four functions from $\cEABC$ to $\cP(N)$}
\label{fig:imset_to_sets}
\end{figure}

Let us fix an arbitrary representation $\repr$ 
and let 
\[
\cEABC^{\repr} = \{u \in \cEABC \mid 0 \leq \exists s \leq \abs{A}-1, 0 \leq \exists t \leq \abs{B}-1 \ s.t. \ u = u_{\repr}(s,t) \in \cEABC^{s,t} \}.
\]
Then there exist $a_{st}, b_{st}$ and $\Gamma_{st} \in \cP(N)$ for every $u_{\repr}(s,t) \in \cEABC^{\repr}$ such that 
$u_{\repr}(s,t) = u_{\tri{a_{st}}{b_{st}}{\Gamma_{st}}}$. 
We define four functions from $\cEABC^{\repr}$ to 
$\mathcal{P}(N)$ as 
$\GammaSetRepr(u_{\repr}(s,t))=\Gamma_{st}$, 
$\aGammaSetRepr(u_{\repr}(s,t))=a_{st}\Gamma_{st}$, 
$\bGammaSetRepr(u_{\repr}(s,t))=b_{st}\Gamma_{st}$ and 
$\abGammaSetRepr(u_{\repr}(s,t))=a_{st}b_{st}\Gamma_{st}$. 
See Figure~\ref{fig:imset_to_sets}. 
We often omit the superscript $\repr$ in $\GammaSetRepr, \aGammaSetRepr, \bGammaSetRepr, \abGammaSetRepr$ and simply write
$\GammaSet, \aGammaSet, \bGammaSet, \abGammaSet$ for notational simplicity.
Then we have the following lemma. 
\begin{lemma}
 \label{lem_tanaken_cor:1}
 For every $v \in \cEABC^{\repr}$, the following equations hold. 
\begin{equation}
\abs{A \cap \GammaSet(v)} + 1 = \abs{A \cap \aGammaSet(v)} = \abs{A \cap \bGammaSet(v)} + 1 = \abs{A \cap \abGammaSet(v)},
 \nonumber 
\end{equation}
\begin{equation}
 \abs{B \cap \GammaSet(v)} + 1 = \abs{B \cap \aGammaSet(v)} + 1 = \abs{B \cap \bGammaSet(v)} = \abs{B \cap \abGammaSet(v)}
  \nonumber 
\end{equation}
\end{lemma}
\begin{proof}
This follows immediately from the definitions of $\GammaSet$, $\aGammaSet$, $\bGammaSet$ and $\abGammaSet$. 
\qed
\end{proof}

Let $\cQ = \cup_{v \in \cEABC^{\repr}}
\{ \GammaSet(v), \aGammaSet(v), \bGammaSet(v), \abGammaSet(v)\}$.
Then we have the following lemma. 
\begin{lemma}
 \label{lem:tanaken_cor:2}
 For any $Q \in \cQ$, let $s = \abs{A \cap Q}$ and $t = \abs{B \cap Q}$. 
 Then each of $\GammaSet^{-1}(Q)$, $\aGammaSet^{-1}(Q)$, $\bGammaSet^{-1}(Q)$ and $\abGammaSet^{-1}(Q)$ 
 has at most one element: 
 $\GammaSet^{-1}(Q) = \emptyset\ or\ \{u(s,t)\}$, 
 $\aGammaSet^{-1}(Q) = \emptyset\ or\ \{u(s-1,t)\}$, 
 $\bGammaSet^{-1}(Q) = \emptyset\ or\ \{u(s,t-1)\}$ and 
 $\abGammaSet^{-1}(Q) = \emptyset\ or\ \{u(s-1,t-1)\}$. 
 Especially,  
 $\GammaSet^{-1}(\GammaSet(v)) = \{v\}$,
 $\aGammaSet^{-1}(\aGammaSet(v)) = \{v\}$,
 $\bGammaSet^{-1}(\bGammaSet(v)) = \{v\}$ and 
 $\abGammaSet^{-1}(\abGammaSet(v)) = \{v\}$ hold. 
\end{lemma}
\begin{proof}
This follows immediately from Proposition~\ref{prop:|A||B|}.
\qed
\end{proof}

For the four ``vertices'' of the representation, we have the following lemma. 
\begin{lemma}
 \label{lem:tanaken_cor:3}
 For the four points $C, AC, BC, ABC \in \cP(N)$, 
 the following four properties hold. 
 \begin{enumerate}
  \item $\GammaSet^{-1}(C)=\{u(0,0)\}$, $\aGammaSet^{-1}(C)=\bGammaSet^{-1}(C)=\abGammaSet^{-1}(C)=\emptyset$
	\label{lem:tanaken_cor:3-1}
  \item $\aGammaSet^{-1}(AC)=\{u(\abs{A}-1,0)\}$, $\GammaSet^{-1}(AC)=\bGammaSet^{-1}(AC)=\abGammaSet^{-1}(AC)=\emptyset$
	\label{lem:tanaken_cor:3-2}
  \item $\bGammaSet^{-1}(BC)=\{u(0, \abs{B}-1)\}$, $\GammaSet^{-1}(BC)=\aGammaSet^{-1}(BC)=\abGammaSet^{-1}(BC)=\emptyset$
	\label{lem:tanaken_cor:3-3}
  \item $\abGammaSet^{-1}(ABC)=\{u(\abs{A}-1, \abs{B}-1)\}$, $\GammaSet^{-1}(ABC)=\aGammaSet^{-1}(ABC)=\bGammaSet^{-1}(ABC)=\emptyset$
	\label{lem:tanaken_cor:3-4}
 \end{enumerate}
\end{lemma}
\begin{proof}
We give the proof for the property in \ref{lem:tanaken_cor:3-2}. 
The proofs for \ref{lem:tanaken_cor:3-1}, \ref{lem:tanaken_cor:3-3} and \ref{lem:tanaken_cor:3-4} are almost the same as in \ref{lem:tanaken_cor:3-2}. 

Let us consider $\GammaSet^{-1}(AC)$. 
Assume that there exists $v \in \cEABC^{\repr}$ such that $v \in \GammaSet^{-1}(AC)$. 
Then, from Lemma~\ref{lem_tanaken_cor:1} and $\abs{A \cap \GammaSet(v)} = \abs{A}$, 
we have $\abs{A \cap \aGammaSet(v)} = \abs{A}+1$. 
This is a contradiction. Therefore such $v \in \cEABC^{\repr}$ does not exist. 

Next, from $B \cap AC = \emptyset$, we have $\bGammaSet^{-1}(AC)=\abGammaSet^{-1}(AC)=\emptyset$. 
Furthermore, because $u_{\tABC}$ takes $-1$ at $AC$ and $\bGammaSet^{-1}(AC)=\emptyset$, 
$\aGammaSet^{-1}(AC)$ must have at least one element. 
Hence, from Lemma~\ref{lem:tanaken_cor:2}, we have $\aGammaSet^{-1}(AC)=\{u(\abs{A}-1,0)\}$. 
\qed
\end{proof}

Let us define the four ``edges'' $\cQ_{ll}, \cQ_{lr}, \cQ_{ul}, \cQ_{ur}$ of the representation as follows: 
\begin{eqnarray}
 \cQ_{ll} & = & \{\aGammaSet(u(s, 0)) \mid 0 \leq s \leq \abs{A}-2 \},
  \nonumber \\
 \cQ_{lr} & = & \{\bGammaSet(u(0, t)) \mid 0 \leq t \leq \abs{B}-2 \},
  \nonumber \\
 \cQ_{ul} & = & \{\aGammaSet(u(\abs{A}-1), t)) \mid 0 \leq t \leq \abs{B}-2 \},
  \nonumber \\
 \cQ_{ur} & = & \{\bGammaSet(u(s, \abs{B}-1))) \mid 0 \leq s \leq \abs{A}-2 \}.
  \nonumber 
\end{eqnarray}
For the four edges, we have the following lemma. 
\begin{lemma}
 \label{lem:tanaken_cor:4}
 For $\cQ_{ll}, \cQ_{lr}, \cQ_{ul}, \cQ_{ur}$, 
 the following four properties hold. 
 \begin{enumerate}
  \item For $u(s,0) \in \cQ_{ll}$, 
	$\GammaSet^{-1}(\aGammaSet(u(s,0))) = \{u(s+1,0)\}$, 
	$\aGammaSet^{-1}(\aGammaSet(u(s,0))) = \{u(s,0)\}$, 
	$\bGammaSet^{-1}(\aGammaSet(u(s,0))) = \emptyset$ and 
	$\abGammaSet^{-1}(\aGammaSet(u(s,0))) = \emptyset$ hold. 
	\label{lem:tanaken_cor:4-1}
  \item For $u(0,t) \in \cQ_{lr}$, 
	$\GammaSet^{-1}(\bGammaSet(u(0,t))) = \{u(0,t+1)\}$, 
	$\aGammaSet^{-1}(\bGammaSet(u(0,t))) = \emptyset$, 
	$\bGammaSet^{-1}(\bGammaSet(u(0,t))) = \{u(0,t)\}$, and 
	$\abGammaSet^{-1}(\bGammaSet(u(0,t))) = \emptyset$ hold. 
	\label{lem:tanaken_cor:4-2}
  \item For $u(s,\abs{B}-1) \in \cQ_{ul}$, 
	$\GammaSet^{-1}(\abGammaSet(u(s,\abs{B}-1))) = \emptyset$, 
	$\aGammaSet^{-1}(\abGammaSet(u(s,\abs{B}-1))) = \emptyset$, 
	$\bGammaSet^{-1}(\abGammaSet(u(s,\abs{B}-1))) = \{u(s+1,\abs{B}-1)\}$ and 
	$\abGammaSet^{-1}(\abGammaSet(u(s,\abs{B}-1))) = \{u(s,\abs{B}-1)\}$ hold. 
	\label{lem:tanaken_cor:4-3}
  \item For $u(\abs{A}-1,t) \in \cQ_{ur}$, 
	$\GammaSet^{-1}(\abGammaSet(u(\abs{A}-1,t))) = \emptyset$, 
	$\aGammaSet^{-1}(\abGammaSet(u(\abs{A}-1,t))) = \{u(\abs{A}-1,t+1)\}$, 
	$\bGammaSet^{-1}(\abGammaSet(u(\abs{A}-1,t))) = \emptyset$ and 
	$\abGammaSet^{-1}(\abGammaSet(u(\abs{A}-1,t))) = \{u(\abs{A}-1,t)\}$ hold. 
	\label{lem:tanaken_cor:4-4}
 \end{enumerate}
\end{lemma}
\begin{proof}
We give the proof for the property in \ref{lem:tanaken_cor:4-1}. 
The proofs for \ref{lem:tanaken_cor:4-2}, \ref{lem:tanaken_cor:4-3} and \ref{lem:tanaken_cor:4-4} are almost the same as in \ref{lem:tanaken_cor:4-1}. 

First, from Lemma~\ref{lem:tanaken_cor:2}, we have $\aGammaSet^{-1}(\aGammaSet(u(s,0))) = \{u(s,0)\}$. 

Next, since $B \cap \aGammaSet(u(s,0)) = \emptyset$, 
we have $\bGammaSet^{-1}(\aGammaSet(u(s,0))) = \abGammaSet^{-1}(\aGammaSet(u(s,0))) = \emptyset$. 

Finally, we consider $\GammaSet^{-1}(\aGammaSet(u(s,0)))$. 
Because $u_{\tABC}$ takes $0$ at $\aGammaSet(u(s,0))$ and $\abGammaSet^{-1}(\aGammaSet(u(s,0))) = \emptyset$, 
$\GammaSet^{-1}(\aGammaSet(u(s,0)))$ must have at least one element to cancel out the value $-1$ of $\aGammaSet^{-1}(\aGammaSet(u(s,0)))$. 
Hence, from Lemma~\ref{lem:tanaken_cor:2}, we have $\GammaSet^{-1}(\aGammaSet(u(s,0))) = \{u(s+1,0)\}$.  
\qed
\end{proof}

\begin{figure}
\begin{center}
\includegraphics[scale=1]{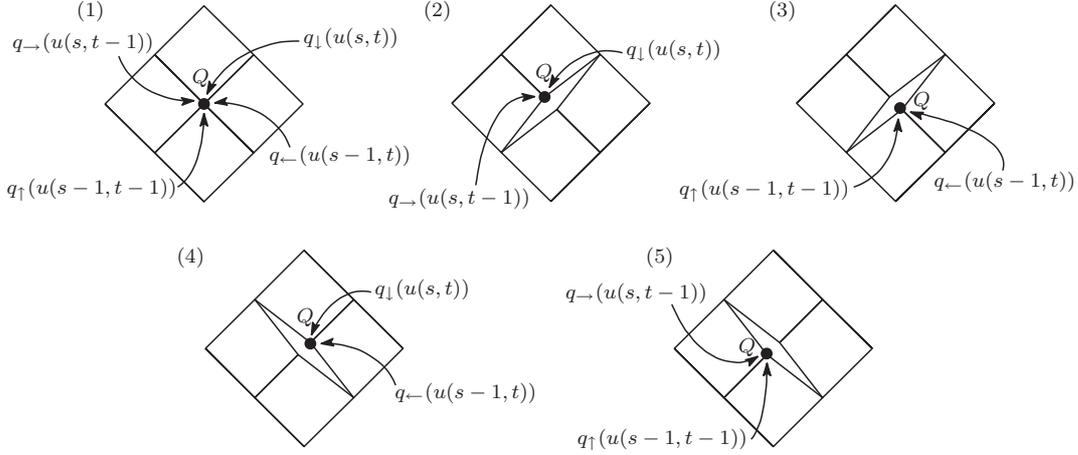}
\end{center}
\caption{Five types of inner points}
\label{fig:five_types_of_rift_patterns}
\end{figure}

For the ``inner points'' of the representation, 
we have the following lemma. 
\begin{lemma}
 \label{lem:tanaken_cor:5}
 Let $\cQ_{i} = \cQ \setminus (\{C, AC, BC, ABC\} \cup \cQ_{ll} \cup \cQ_{lr} \cup \cQ_{ul} \cup \cQ_{ur})$,
 and $s = \abs{A \cap Q}$, $t = \abs{B \cap Q}$ for $Q \in \cQ_{i}$. 
 For each $Q \in \cQ_{i}$, 
 one of the following five properties holds. 
 \begin{enumerate}
  \item \label{lem:tanaken_cor:5-1}
	$\GammaSet^{-1}(Q) = \{u(s,t)\}$, $\aGammaSet^{-1}(Q) = \{u(s-1,t)\}$, $\bGammaSet^{-1}(Q) = \{u(s,t-1)\}$, $\abGammaSet^{-1}(Q) = \{u(s-1,t-1)\}$
	((1) in Figure~\ref{fig:five_types_of_rift_patterns})
  \item \label{lem:tanaken_cor:5-2}
	$\GammaSet^{-1}(Q) = \{u(s,t)\}$, $\aGammaSet^{-1}(Q) = \emptyset$, $\bGammaSet^{-1}(Q) = \{u(s,t-1)\}$, $\abGammaSet^{-1}(Q) = \emptyset$
	((2) in Figure~\ref{fig:five_types_of_rift_patterns})
  \item \label{lem:tanaken_cor:5-3}
	$\GammaSet^{-1}(Q) = \emptyset$, $\aGammaSet^{-1}(Q) = \{u(s-1,t)\}$, $\bGammaSet^{-1}(Q) = \emptyset$, $\abGammaSet^{-1}(Q) = \{u(s-1,t-1)\}$
	((3) in Figure~\ref{fig:five_types_of_rift_patterns})
  \item \label{lem:tanaken_cor:5-4}
	$\GammaSet^{-1}(Q) = \{u(s,t)\}$, $\aGammaSet^{-1}(Q) = \{u(s-1,t)\}$, $\bGammaSet^{-1}(Q) = \emptyset$, $\abGammaSet^{-1}(Q) = \emptyset$
	((4) in Figure~\ref{fig:five_types_of_rift_patterns})
  \item \label{lem:tanaken_cor:5-5}
	$\GammaSet^{-1}(Q) = \emptyset$, $\aGammaSet^{-1}(Q) = \emptyset$, $\bGammaSet^{-1}(Q) = \{u(s,t-1)\}$, $\abGammaSet^{-1}(Q) = \{u(s-1,t-1)\}$
	((5) in Figure~\ref{fig:five_types_of_rift_patterns})
 \end{enumerate}
\end{lemma}
\begin{proof}
Let $Q \in \cQ_{i}$. 
First, we show that at least one of $\GammaSet^{-1}(Q)$ and $\abGammaSet^{-1}(Q)$ is not empty. 
If $\GammaSet^{-1}(Q)$ and $\abGammaSet^{-1}(Q)$ are both empty, 
then $\aGammaSet^{-1}(Q)$ and $\bGammaSet^{-1}(Q)$ must be empty 
because 
i) the value at $Q \in \cQ_{i}$ is zero, 
ii) $\GammaSet^{-1}(Q)$ and $\abGammaSet^{-1}(Q)$ have the value $1$ by (\ref{eq:semi-elementary})
and iii) $\aGammaSet^{-1}(Q)$ and $\bGammaSet^{-1}(Q)$ have the value $-1$ by (\ref{eq:semi-elementary}). 
This contradicts $Q \in \cQ_{i}$. 
Then, from Lemma~\ref{lem:tanaken_cor:2}, 
the possible patterns are classified as follows. 
\renewcommand{\labelenumi}{(\Roman{enumi})}
\begin{enumerate}
 \item We consider the case where both $\GammaSet^{-1}(Q)$ and $\abGammaSet^{-1}(Q)$ are non-empty. 
       In this case, 
       $\aGammaSet^{-1}(Q)$ and $\bGammaSet^{-1}(Q)$ must be non-empty 
       in order to cancel out the value $2$ by $\GammaSet^{-1}(Q)$ and $\abGammaSet^{-1}(Q)$ at $Q$. 
       This is the case of \ref{lem:tanaken_cor:5-1} in Lemma~\ref{lem:tanaken_cor:5}. 
 \item We consider the case where $\GammaSet^{-1}(Q)$ is non-empty and $\abGammaSet^{-1}(Q)$ is empty. 
       In this case, one of $\aGammaSet^{-1}(Q)$ and $\bGammaSet^{-1}(Q)$ must be non-empty
       and another must be empty in order to cancel out the value $1$ at $Q$. 
       If $\aGammaSet^{-1}(Q)$ is non-empty and $\bGammaSet^{-1}(Q)$ is empty, 
       then this is the case of \ref{lem:tanaken_cor:5-4} in Lemma~\ref{lem:tanaken_cor:5}. 
       If $\aGammaSet^{-1}(Q)$ is empty and $\bGammaSet^{-1}(Q)$ is non-empty, 
       then this is the case of \ref{lem:tanaken_cor:5-2} in Lemma~\ref{lem:tanaken_cor:5}. 
 \item We consider the case where $\GammaSet^{-1}(Q)$ is empty and $\abGammaSet^{-1}(Q)$ is non-empty. 
       In this case, one of $\aGammaSet^{-1}(Q)$ and $\bGammaSet^{-1}(Q)$ must be non-empty
       and another must be empty in order to cancel out the value $1$ at $Q$. 
       If $\aGammaSet^{-1}(Q)$ is non-empty and $\bGammaSet^{-1}(Q)$ is empty, 
       then this is the case of \ref{lem:tanaken_cor:5-3} in Lemma~\ref{lem:tanaken_cor:5}. 
       If $\aGammaSet^{-1}(Q)$ is empty and $\bGammaSet^{-1}(Q)$ is non-empty, 
       then this is the case of \ref{lem:tanaken_cor:5-5} in Lemma~\ref{lem:tanaken_cor:5}.
\end{enumerate}
\renewcommand{\labelenumi}{\arabic{enumi}}
\qed
\end{proof}

In order to prove Proposition~\ref{prop:sufficient_condition_for_sigma_decomposability}, we use the following lemma. 
\begin{lemma}
 \label{lemma:separable-char}
 Let $1\le s_0\le |A|-1$.
 Then (\ref{eq:ust}) is $(s_0,*)$-separable if and only if
 $\GammaSetSei(u(s_0,t))=\bGammaSetSei(u(s_0,t-1))$ for all $1\le t\le |B|-1$.
 Similarly, for $1\le t_0\le |B|-1$, (\ref{eq:ust}) is $(*,t_0)$-separable if and only if
 $\GammaSetSei(u(s,t_0))=\aGammaSetSei(u(s-1,t_0))$ for all $1\le s\le |B|-1$.
\end{lemma}
\begin{proof}
 This follows immediately from Lemma~\ref{lem:tanaken_cor:4} and \ref{lem:tanaken_cor:5}. 
\qed
\end{proof}

\begin{flushleft}
 {\bf Proof of Proposition~\ref{prop:sufficient_condition_for_sigma_decomposability}}
\end{flushleft}
\begin{proof}
 We prove the proposition by induction.
 If $|A|=|B|=1$, then (\ref{eq:ust}) is $\sigma$-decomposable
 from the definition.
 If $|A|=1$ and $|B|\ge 2$,
 then (\ref{eq:ust}) is $(*,t_0)$-separable for any $1\leq t_0\leq |B|-1$
 and therefore $\sigma$-decomposable by induction.
 The case $|B|=1$ is similarly proved.
 Now assume $|A|=2$ and $|B|\ge 2$.
 If (\ref{eq:ust}) is $(1,*)$-separable,
 then the two representations in (\ref{eq:t-sep}) after separation
 are $\sigma$-decomposable.
 Hence (\ref{eq:ust}) is $\sigma$-decomposable.
 Now assume (\ref{eq:ust}) is not $(1,*)$-separable.
 Then, by Lemma~\ref{lemma:separable-char}, $\GammaSetSei(u(1,t))\neq \bGammaSetSei(u(1,t-1))$ for any $1\le t\le |B|-1$.
 Therefore $\GammaSetSei(u(1,t))=\aGammaSetSei(u(0,t))$ for any $1\le t\le |B|-1$.
 This means that, by Lemma~\ref{lemma:separable-char} again, (\ref{eq:ust})
 is $(*,t)$-separable for any $1\le t\le |B|-1$.
 By induction, we obtain $\sigma$-decomposability of (\ref{eq:ust}).
 The case $|A|\ge 2$ and $|B|=2$ is similarly proved.
\qed
\end{proof}

Next, we consider the rigorous definition of the ``rifts'' of $u_{\tABC}$. 
Let $\cQ_{1}$ be 
\begin{eqnarray}
 \cQ_{1} & = & 
   \{Q \in \cQ \mid 1 \leq \exists s \leq \abs{A}-1, 1 \leq \exists t \leq \abs{B}-1,\, s.t. 
  \nonumber \\ & &  \qquad \qquad \qquad \qquad
  Q=\GammaSet(u(s, t))=\abGammaSet(u(s-1, t-1)) \}. 
 \nonumber 
\end{eqnarray}
Note that $Q=\GammaSet(u(s, t))=\abGammaSet(u(s-1, t-1))$ is equivalent to 
$Q=\GammaSet(u(s, t))=\aGammaSet(u(s-1, t))=\bGammaSet(u(s, t-1))=\abGammaSet(u(s-1, t-1))$. 
Hence each element of $\cQ_{1}$ satisfies \ref{lem:tanaken_cor:5-1} of Lemma~\ref{lem:tanaken_cor:5}. 
Let $\cQ'$ be  
\begin{eqnarray}
 \cQ' 
  & = & 
  \{C, AC, BC, ABC\} \cup \cQ_{ll} \cup \cQ_{lr} \cup \cQ_{ul} \cup \cQ_{ur}
  \cup 
  \cQ_{1}. 
  \nonumber 
\end{eqnarray}
Furthermore, let $\cQ'' = \cQ \setminus \cQ'$. 
Then we have 
\begin{eqnarray}
  \cQ'' & = & \{\GammaSet(u(s, t)), \abGammaSet(u(s-1, t-1)) \mid 1 \leq \exists s \leq \abs{A}-1, 1 \leq \exists t \leq \abs{B}-1,\;  s.t. \quad 
  \nonumber \\ 
  & & \qquad  \qquad \qquad \qquad \qquad \qquad \qquad \qquad \qquad
   \GammaSet(u(s, t)) \neq \abGammaSet(u(s-1, t-1)) \}. 
   \nonumber 
\end{eqnarray}
Each element of $\cQ''$ satisfies 
one of \ref{lem:tanaken_cor:5-2}, \ref{lem:tanaken_cor:5-3}, \ref{lem:tanaken_cor:5-4} and  \ref{lem:tanaken_cor:5-5} 
of Lemma~\ref{lem:tanaken_cor:5}.

Let $\cardA(Q) = \abs{A \cap Q}$ and $\cardB(Q) = \abs{B \cap Q}$ for $Q \in \cQ$.
Now we define {\srift s} and {\brift s} for a representation of $u_{\tABC}$. 
\begin{definition}
 \label{def:rift}
For any $Q \in \cQ''$ such that $Q = \GammaSet(u(s,t)) \neq \aGammaSet(u(s-1,t))$, 
let $t_{L}=t_{L}(Q)$ be the maximum number $t'$ such that $\GammaSet(u(s,t')) = \aGammaSet(u(s-1,t'))$ and $t' \leq t$. 
Furthermore, let $t_{U}=t_{U}(Q)$ be the minimum number $t'$ such that $\bGammaSet(u(s,t'-1)) = \abGammaSet(u(s-1,t'-1))$ and $t' \geq t$. 
Then we define $\rift_{s}(s;t_{L},t_{U}) = \{Q \in \cQ \mid \cardA(Q)=s, t_{L} \leq \cardB(Q) \leq t_{U}\}$
and call it a {\em \srift}. 
See the left-hand side of Figure~\ref{fig:rift}. 

We also define another type of rift. 
For any $Q \in \cQ''$ such that $Q = \GammaSet(u(s,t)) \neq \bGammaSet(u(s-1,t))$, 
let $s_{L}=s_{L}(Q)$ be the maximum number $s'$ such that $\GammaSet(u(s',t)) = \bGammaSet(u(s',t-1))$ and $s' \leq s$. 
Furthermore, let $s_{U}=s_{U}(Q)$ be the minimum number $s'$ such that $\aGammaSet(u(s'-1,t)) = \abGammaSet(u(s'-1,t-1))$ and $t' \geq t$. 
Then we define $\rift_{b}(s_{L},s_{U};t) = \{Q \in \cQ \mid s_{L} \leq \cardA(Q) \leq s_{U}, \cardB(Q)=t\}$ 
and call it a {\em \brift}. 
See the right-hand side of Figure~\ref{fig:rift}. 

A {\em rift} is any \srift\ or \brift,
and {\em a rift pattern} is a set of rifts which do not intersect each other.
\begin{figure}
\begin{center}
\includegraphics[scale=1]{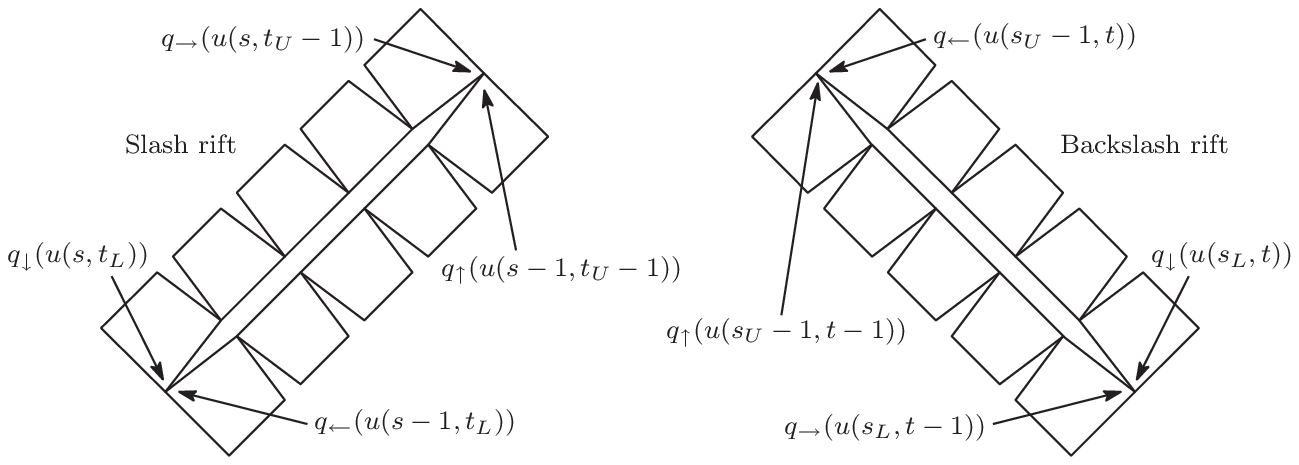}
\end{center}
\caption{Definitions of \srift\ and \brift}
\label{fig:rift}
\end{figure}
\end{definition}
Note that the existence of the maximizer and minimizer in Definition \ref{def:rift} is assured by Lemma~\ref{lem:tanaken_cor:4}. 
We denote the set of all of the \srift s by $\cR_{s} = \{ \rift_{s}(\cardA(Q);t_{L}(Q),t_{U}(Q)) \mid Q \in \cQ'' \}$ 
and the set of all of the \brift s by $\cR_{b} = \{ \rift_{b}(s_{L}(Q),s_{U}(Q);t(Q)) \mid Q \in \cQ'' \}$. 
Furthermore, note that rift occurs only at $(\abs{A}-1)(\abs{B}-1)$ pairs of $(\cardA(Q),\cardB(Q))$ for $Q \in \cQ''$. 

\begin{flushleft}
{\bf Proof of Proposition~\ref{prop:deg_of_free} }
\end{flushleft}
\begin{proof}
The degree of freedom for any rift with length $l$ is given by 
subtracting the number of cases where the 
rift is collapsed in the middle i.e., 
$\GammaSet(Q) = \abGammaSet(Q)$ holds for a $Q$ in the inner points of the rift, 
by permutations, 
from the number of permutations $l!$. 
Therefore, we obtain (\ref{eq:deg_of_free}). 
\qed
\end{proof}
\begin{flushleft}
{\bf Proof of Proposition~\ref{prop:num_of_rep_for_a2} }
\end{flushleft}
\begin{proof}
 We will use induction on $m$. 
 First, it is easily seen that $r_{2}(1) = 1$. 
 Now we assume that~(\ref{eq:num_of_rep_for_a2}) is true for $m-1$. 
 If a  \brift\ first appears in $\GammaSet(u(1,k))$, i.e., there is no \brift\ in $\GammaSet(u(1,0)), \dots, \GammaSet(u(1,k-1))$, 
 then such a representation can be split in two representations: 
 one is represented as $\sum_{i=0}^{2}\sum_{j=0}^{k-1}u(i,l)$ 
 and another is represented as $\sum_{i=0}^{2}\sum_{j=k}^{m}u(i,l)$. 
 In this case, 
 since the former has $k!$ representations 
 and the latter has $r_{2}(m-k)$ representations from the assumption, 
 there are $k!r_{2}(m-k)$ representations. 
 Therefore, by summing the number of representations for $k \in [1, m-1]$ 
 and $m!$ which is the number of representations with no \brift, 
 we obtain (\ref{eq:num_of_rep_for_a2}). 
\qed
\end{proof}

The following lemma implies that there exists at least one rift which can be eliminated for any representation of $u_{\tABC}$. 
\begin{proposition}
 \label{prop:tanaken_prop:2}
 Assume that $\cR_{s}$ or $\cR_{b}$ is non-empty. 
 Let us denote by $\cR_{s}'$
 the subset of all the 
 $\rift_{s}(s';t_{L}',t_{U}') \in \cR_{s}$
 such that, for any $\rift_{s}(s'';t_{L}'',t_{U}'') \in \cR_{s}$, 
 neither $s'' < s'$ nor $(t_{L}'',t_{U}'') \cap (t_{L}',t_{U}') \neq \emptyset$ holds
 where $(t_{L}',t_{U}')$ means open interval. 
 Let us denote by $\cR_{b}'$
 the subset of all the 
 $\rift_{b}(s_{L}',s_{U}';t') \in \cR_{b}$
 such that, for any $\rift_{b}(s_{L}'',s_{U}'';t'') \in \cR_{b}$, 
 neither $t'' < t'$ nor $(s_{L}'',s_{U}'') \cap (s_{L}',s_{U}') \neq \emptyset$ holds. 
 Then at least one of the following two properties holds. 
 \begin{enumerate}
  \item \label{prop:tanaken_prop:2-1}
	There exists $\rift_{s}(s;t_{L},t_{U}) \in \cR_{s}'$ 
	such that, for any $\rift_{b}(s_{L},s_{U};t) \in \cR_{b}$, either $t \notin [t_{L}+1, t_{U}-1]$ or $s \leq s_{L}$ holds. 
  \item \label{prop:tanaken_prop:2-2}
	There exists $\rift_{b}(s_{L},s_{U};t) \in \cR_{b}'$
	such that, for any $\rift_{s}(s;t_{L},t_{U}) \in \cR_{s}$, either $s \notin [s_{L}+1,s_{U}-1]$ or $t \leq t_{L}$ holds. 
 \end{enumerate}
\end{proposition}
\begin{proof}
If $\cR_{s}$ is non-empty, 
then $\rift_{s}(\cardA(Q);t_{L}(Q),t_{U}(Q))$ with the smallest $s$ in $\cR_{s}$ 
is clearly one of the element of $\cR_{s}'$. 
This means that $\cR_{s}'$ is non-empty if $\cR_{s}$ is non-empty. 
The same thing holds for $\cR_{b}$ and $\cR_{b}'$. 

If either $\cR_{s}'$ or $\cR_{b}$ is empty, 
then either \ref{prop:tanaken_prop:2-1} or \ref{prop:tanaken_prop:2-2} 
in the statement of Proposition~\ref{prop:tanaken_prop:2-1} holds. 

In the following, we consider the case where neither $\cR_{s}'$ nor $\cR_{b}$ is empty. 
We show by contradiction that either \ref{prop:tanaken_prop:2-1} or 
\ref{prop:tanaken_prop:2-2} holds.  Suppose that neither \ref{prop:tanaken_prop:2-1}
nor \ref{prop:tanaken_prop:2-2} holds. 

Now we choose one element from $\cR_{s}'$ and denote it by $\rift_{s}(s^{(1)};t_{L}^{(1)},t_{U}^{(1)})$. 
By assumption that {\ref{prop:tanaken_prop:2-1}} does not hold, 
we can choose $\rift_{b}(s_{L},s_{U};t) \in \cR_{b}$ 
such that $t \in [t_{L}^{(1)}+1,t_{U}^{(1)}-1]$ and $s^{(1)} > s_{L}$ hold for $\rift_{s}(s^{(1)};t_{L}^{(1)},t_{U}^{(1)})$. 

If $\rift_{b}(s_{L},s_{U};t) \in \cR_{b}'$ holds, then we denote it by $\rift_{b}(s_{L}^{(1)},s_{U}^{(1)};t^{(1)})$. 
If not, we choose $\rift_{b}(s_{L}^{(1)},s_{U}^{(1)};t^{(1)}) \in \cR_{b}'$ 
such that $t^{(1)} < t$ and $(s_{L}^{(1)},s_{U}^{(1)}) \cap (s_{L}, s_{U}) \neq \emptyset$ hold. 
In both cases, we have $t^{(1)} \leq t_{U}^{(1)} - 1$. 

Next, 
by assumption that {\ref{prop:tanaken_prop:2-1}} does not hold, 
we can choose one element from $\rift_{s}(s;t_{L},t_{U}) \in \cR_{s}$ 
such that $s \in [s_{L}^{(1)}+1,s_{U}^{(1)}-1]$ and $t^{(1)} > t_{L}$ hold 
for $\rift_{b}(s_{L}^{(1)},s_{U}^{(1)};t^{(1)})$. 
In this case, we have $t^{(1)} \geq t_{U}$. 
This can be shown as follows. 
If $t^{(1)} \geq t_{U}$ does not hold, i.e., $t_{L} < t^{(1)} < t_{U}$, 
then $Q=\GammaSet(u(s,t^{(1)}))$ belongs to 
both $\rift_{s}(s;t_{L},t_{U})$ and $\rift_{b}(s_{L}^{(1)},s_{U}^{(1)};t^{(1)})$. 
This contradicts Lemma~\ref{lem:tanaken_cor:5} and hence we have $t^{(1)} \geq t_{U}$. 

If $\rift_{s}(s;t_{L},t_{U}) \in \cR_{s}'$ holds, then we denote it by $\rift_{s}(s^{(2)};t_{L}^{(2)},t_{U}^{(2)}) \in \cR_{s}'$. 
If not, we choose $\rift_{s}(s^{(2)};t_{L}^{(2)},t_{U}^{(2)}) \in \cR_{s}'$ 
such that $s^{(2)} < s$ and $(t_{L}^{(2)},t_{U}^{(2)}) \cap (t_{L}, t_{U}) \neq \emptyset$ hold. 
In both cases, we have $s^{(2)} \leq s_{U}^{(1)}-1$. 
Furthermore, by considering the relation $t_{L}^{(2)} \leq t_{U} \leq t^{(1)} \leq t_{U}^{(1)} - 1$
and the fact that 
both $\rift_{s}(s^{(1)};t_{L}^{(1)},t_{U}^{(1)})$ and $\rift_{s}(s^{(2)};t_{L}^{(2)},t_{U}^{(2)})$ 
are the elements of $\cR_{s}'$, 
we have $(t_{L}^{(2)},t_{U}^{(2)}) \cap (t_{L}^{(1)},t_{U}^{(1)})=\emptyset$. 
This means that $t_{U}^{(2)} \leq t_{L}^{(1)}$ holds. 
By applying much the same way as above to $\rift_{s}(s^{(2)};t_{L}^{(2)},t_{U}^{(2)})$, 
we can choose 
$\rift_{b}(s_{L}^{(2)},s_{U}^{(2)};t^{(2)}) \in \cR_{b}'$
such that $s_{U}^{(2)} \leq s_{L}^{(1)}$ holds. 

Therefore, for 
$\rift_{s}(s^{(1)};t_{L}^{(1)},t_{U}^{(1)})$, $\rift_{s}(s^{(2)};t_{L}^{(2)},t_{U}^{(2)})$
$\rift_{b}(s_{L}^{(1)},s_{U}^{(1)};t^{(1)})$ and $\rift_{s}(s^{(2)};t_{L}^{(2)},t_{U}^{(2)})$, 
we obtain the following inequalities: 
\begin{equation}
 s_{L}^{(2)} < s_{U}^{(2)} \leq s_{L}^{(1)} < s_{U}^{(1)}
  ,\ 
 t_{L}^{(2)} < t_{U}^{(2)} \leq t_{L}^{(1)} < t_{U}^{(1)}
 \label{eq:st_loop}
\end{equation}

If neither \ref{prop:tanaken_prop:2-1} nor \ref{prop:tanaken_prop:2-2} holds, 
we can repeat the above procedure indefinitely. 
However, from~(\ref{eq:st_loop}), 
$s$ and $t$ decrease in each step. 
This means that the same rift never appear more than once in the procedure. 
The procedure stops in finite times and this is contradiction.  
Therefore either \ref{prop:tanaken_prop:2-1} or \ref{prop:tanaken_prop:2-2} 
in the statement of Proposition~\ref{prop:tanaken_prop:2} holds. 
\qed
\end{proof}

Now we give  the proof of Theorem~\ref{thm:main}. 
\begin{flushleft}
{\bf Proof of Theorem~\ref{thm:main} }
\end{flushleft}
\begin{proof}
It suffices to show that 
 any representation of $u_{\tABC}$ and the standard representation in (\ref{eq:standard_represetation})
 are connected by two-by-two basic relations. 
 Remember, from Proposition~\ref{prop:kashimura} in Section~\ref{sec:main results}, that
 the representations without rifts are connected by two-by-two basic relations. 
 Therefore, from Proposition~\ref{prop:tanaken_prop:2}, 
 it suffices to show that 
 we can eliminate the rift 
 which satisfies either \ref{prop:tanaken_prop:2-1} or \ref{prop:tanaken_prop:2-2} 
 in the statement of Proposition~\ref{prop:tanaken_prop:2}. 
 We prove this in the case where a rift satisfies
 \ref{prop:tanaken_prop:2-2} of Proposition~\ref{prop:tanaken_prop:2}. 
 The proof for the case of \ref{prop:tanaken_prop:2-1} 
 is similar to the proof of the case of \ref{prop:tanaken_prop:2-2} and will be omitted. 

 Let us consider a representation: 
 \begin{equation}
  \label{eq:init_representation}
  u_{\tABC}=\sum_{s=0}^{|A|-1} \sum_{t=0}^{|B|-1} u_{\repr}(s,t). 
 \end{equation}
 In the rest of the proof, we omit the subscript $\repr$ in $u_{\repr}(s,t)$.
 Assume that a rift in this representation satisfies
 \ref{prop:tanaken_prop:2-2} of Proposition~\ref{prop:tanaken_prop:2}. 
 Then there exists $\rift_{b}(s_{L},s_{U};t) \in \cR_{b}$ 
 such that either $s \notin [s_{L}+1,s_{U}-1]$ or $t \leq t_{L}$ 
 holds for any $\rift_{s}(s;t_{L},t_{U}) \in \cR_{s}$. 
 We show that we can eliminate the \brift\ $\rift_{b}(s_{L},s_{U};t)$
 by applying two-by-two basic relations. 

 Let 
 $\Gamma' = \GammaSet(u(s_{L},0))$, 
 $B' = B \cap \bGammaSet(u(s_{L},t-1))$ and  
 $A' = A \cap \{\aGammaSet(u(s_{U}-1,0)) \setminus \Gamma'\}$. 
 Then we obtain the following representation of $u_{\tri{A'}{B'}{\Gamma'}}$:
 \begin{equation}
  \label{eq:rep_for_chip}
   u_{\tri{A'}{B'}{\Gamma'}} = \sum_{s' = s_{L}}^{s_{U}-1} \sum_{t'=0}^{t-1} u(s',t'). 
 \end{equation} 
 From the assumption, i.e., \ref{prop:tanaken_prop:2-2} of Proposition~\ref{prop:tanaken_prop:2}, 
 it is easily seen that the representation in (\ref{eq:rep_for_chip}) has no rift. 
 Let
 ${a}^{(s')} = \aGammaSet(u(s'+s_{L}-1,0)) \setminus \GammaSet(u(s'+s_{L}-1,0))$ 
 for $1 \leq s' \leq s_{U}-s_{L}$. 
 Note that the sequence ${a}^{(1)},\dots,{a}^{(s_{U}-s_{L})}$ 
 is ordered with respect to the upper side of the rift. 
 In the same way, let
 ${b}^{(t')} = \bGammaSet(u(s_{L}, t'-1)) \setminus \GammaSet(u(s_{L}, t'-1))$
 for $1 \leq t' \leq t$. 
 Let us consider the following representation without rifts: 
 \begin{equation}
  \label{eq:rep_for_chip_after_transformation}
   u_{\tri{A'}{B'}{\Gamma'}}
   = 
   \sum_{s' = 1}^{s_{U}-s_{L}} \sum_{t'=1}^{t} 
    u_{\tri{a^{(s')}}{b^{(t')}}{a^{(1)} \dots a^{(s'-1)} b^{(1)} \dots b^{(t'-1)}\Gamma'}}.
 \end{equation}
 
 Remember that the terms in the representation of (\ref{eq:rep_for_chip}) 
 also appear in the representation of (\ref{eq:init_representation}). 
 Hence, in the representation of (\ref{eq:init_representation}), 
 if we replace the terms used in the representation of (\ref{eq:rep_for_chip}) 
 with the terms used in the representation of (\ref{eq:rep_for_chip_after_transformation}), 
 then we can eliminate the \brift\ $\rift_{b}(s_{L},s_{U};t)$ 
 in the representation of (\ref{eq:init_representation})
 because of the definition of ${a}^{(s')}$ and ${b}^{(t')}$. 
 Furthermore, from Proposition~\ref{prop:kashimura}, 
 this replacement is realized by applying two-by-two basic relations. 
 Therefore we can eliminate the \brift\ $\rift_{b}(s_{L},s_{U};t)$ by applying two-by-two basic relations. 

 This completes the proof of Theorem~\ref{thm:main}.
\qed
\end{proof}

\bibliography{imset}
\bibliographystyle{plain}
\end{document}